\documentclass[12pt]{article}
\usepackage{timet,color}
\usepackage[urlcolor=blue,citecolor=black,linkcolor=black]{hyperref}

\usepackage{natbib}
\usepackage{graphicx} 
\usepackage{amsmath} 
\usepackage{amssymb}
\usepackage{booktabs}
\usepackage{multirow}

\title{Fast parallel transient electromagnetic modeling using a~uniform-in-time approximation to the exponential}

\author{Ralph-Uwe Börner\footnote{Institute of Geophysics and Geoinformatics, TU Bergakademie Freiberg, Germany}  \and Stefan Güttel\footnote{Department of Mathematics, The University of Manchester, United Kingdom}}

\date{April 2025}

\newcommand{\mA}{\mathbf A}
\newcommand{\mI}{\mathbf I}
\newcommand{\mJ}{\mathbf J}
\newcommand{\mK}{\mathbf K}
\newcommand{\mM}{\mathbf M}

\newcommand{\mQ}{\mathbf Q}
\newcommand{\mS}{\mathbf S}
\newcommand{\mR}{\mathbf R}
\newcommand{\vb}{\mathbf b}
\newcommand{\vf}{\mathbf f}
\newcommand{\vg}{\mathbf g}
\newcommand{\vu}{\mathbf u}
\newcommand{\ve}{\mathbf e}
\newcommand{\vm}{\mathbf m}
\newcommand{\vr}{\mathbf r}
\newcommand{\vv}{\mathbf v}
\newcommand{\vx}{\mathbf x}
\newcommand{\vy}{\mathbf y}

\usepackage[font=small,labelfont=bf]{caption}

\begin{document}

\maketitle

\begin{abstract}
    A new approach for the parallel forward modeling of transient electromagnetic (TEM) fields is presented. It is based on a family of uniform-in-time rational approximants to the matrix exponential that share a common denominator independent of the evaluation time points. The partial fraction decomposition of this family is exploited to devise a fast solver with high parallel efficiency. The number of shifted linear systems that need to be solved in parallel does not depend on the number of required time channels nor the spatial discretization. We also argue that similar parallel efficiency gains can be expected when solving the inverse TEM problem.
\end{abstract}


\sloppy
\section{Introduction}

The numerical simulation of transient electromagnetic (TEM) fields in three-dimensional conductivity models plays a central role in the interpretation of TEM data. In particular, the accurate and efficient solution of the forward problem is a critical component of inversion workflows.

Interest in numerical approaches to the TEM forward problem dates back to the 1980s, with early studies such as those by \citet{oristaglio1982diffusion,oristaglio1984diffusion} who employed finite difference grids to model TEM field propagation in two-dimensional conductivity structures.
Among the earliest and most widely adopted numerical schemes is the finite-difference time-domain (FDTD) method, which solves Maxwell’s equations using explicit time integration, notably the Du Fort--Frankel scheme proposed by \citet{dufort1953stability}. This approach has been extended to three-dimensional settings; see, e.g., \citet{wang1993finite}. The key advantage of these explicit schemes lies in their computational simplicity, as each time step involves only inexpensive matrix-vector products. However, their applicability is often constrained by stringent stability conditions which necessitate extremely small time steps---particularly when modeling fine spatial resolutions or low-conductivity media---leading to high overall computational cost.

Implicit time-integration schemes offer a viable alternative. Although each time step involves the solution of a linear system---thus incurring higher per-step computational cost---implicit methods allow for significantly larger time steps, often resulting in faster simulations overall. \citet{haber20023d} demonstrated the application of an implicit scheme for simulating the 3-D transient electromagnetic response of conductive earth models.

Implicit methods are often favored in inversion contexts because of their flexibility and relatively straightforward implementation. Their computational efficiency can be further enhanced when matrix factorizations (e.g., LU or Cholesky) are precomputed and reused across time steps, reducing the cost of each step to a pair of inexpensive forward and backward substitutions~\citep{um20103d}. 
To further alleviate the computational burden of solving large sparse systems, strategies such as the direct splitting (DS) method have been proposed. This approach transforms large sparse matrices into a sequence of smaller, diagonally dominant tridiagonal systems, thereby improving numerical efficiency~\citep{liu2022}.

An alternative to traditional time-stepping approaches for solving the TEM forward problem is to compute the action of the matrix exponential on a source vector directly, thereby avoiding time discretization altogether.
Krylov subspace methods have proven particularly effective for this purpose.
These methods approximate the solution of the system governed by a parabolic-type partial differential equation by projecting the exponential operator onto a Krylov space and computing the exponential of a much smaller matrix resulting from this projection.
The class of exponential integrators using rational Krylov spaces is particularly attractive for TEM applications because it permits large effective time steps without compromising stability or accuracy, and avoids the numerical stiffness that often affects conventional time-stepping schemes. 
Pioneering work by \citet{druskin2009solution,borner2015three} demonstrated the effectiveness of rational Krylov space methods in geo-electromagnetics, showing that only a modest number of Krylov vectors is typically needed to capture the essential dynamics, particularly at late times when the field responses are smoother. While there are ways to parallelize rational Krylov methods, numerical stability issues may constrain the achievable level of parallel efficiency; see~\cite{BG17a}. 
In recent work, \citet{gao2024} extended the rational Krylov space framework to accommodate arbitrary current waveforms, thereby broadening its applicability to more general transient excitation scenarios. 

An alternative approach to time-domain modeling of transient electromagnetic (TEM) responses is based on frequency-domain solutions synthesized via inverse Fourier (or Laplace) transform techniques.
In this framework, the time-domain response is reconstructed by superposing solutions obtained at a discrete set of frequencies.
This approach allows for the use of mature frequency-domain solvers, which benefit from established numerical techniques. However, accurately capturing the transient response requires a sufficiently dense sampling of frequencies, particularly to resolve sharp transients at early times.
While frequency-domain solvers are naturally parallelizable across frequencies and avoid the time-stepping constraints of explicit schemes, their efficiency depends critically on the choice of sampling frequencies  and the stability of the underlying frequency-domain solver~\citep{rochlitz2021}.

The remainder of this paper is organized as follows. 
In Section~\ref{sec:formulation}, we present the mathematical formulation of the transient electromagnetic problem and introduce the variational framework used for spatial discretization. Section~\ref{sec:numerical_method} details the proposed rational approximation strategies.
Section~\ref{sec:inversion} discusses the implications of the proposed method for inverse modeling, with a focus on its potential for computational acceleration.
In Section~\ref{sec:experiments} we demonstrate the accuracy and efficiency of the method through numerical experiments on a representative test model.
Finally, in Section~\ref{sec:conclusion}, we summarize the main findings and discuss extensions of the approach.

\section{Mathematical model and discretization}
\label{sec:formulation}

\subsection{Governing equations}

We recall the governing equations of electromagnetic induction. 
In the quasi-static approximation and after eliminating the magnetic field, Maxwell's equations for the electric field $\ve = \ve(\vr, t)$ in the time domain read
\begin{equation}
    \label{eq:pde}
    \nabla \times \mu^{-1} \nabla \times \ve + \sigma \partial_t \ve = -\partial_t \mathbf{j}^e, \quad t \in \mathbb{R},
\end{equation}
where the spatial coordinate $\vr$ varies in the computational domain $\Omega \subset \mathbb{R}^3$.
The magnetic permeability is $\mu = \mu_0 = 4 \pi \times 10 ^{-7}$ Vs/(Am), and the electrical conductivity $\sigma = \sigma(\vr)$ is a function defined on $\Omega$.
On the boundary $\partial\Omega$ of the domain $\Omega$ we impose a boundary condition of the form
\begin{equation}
    \mathbf n \times \ve = \mathbf{0}.
\end{equation}
As source term we consider a current density $\mathbf{j}^e$ resulting from a transmitter with a stationary current that is shut off at time $t = 0$,  resulting in 
\begin{equation}
    \mathbf{j}^e(\vr, t) = \mathbf q(\vr)H(-t),
\end{equation}
where $\mathbf q(\vr)$ denotes the spatial trace of the transmitter, and $H(\cdot)$ denotes the Heaviside step function.
In summary, the initial-value problem to be solved is  
\begin{subequations}
    \begin{align}
    \nabla \times \mu^{-1} \nabla \times \ve + \sigma \partial_t \ve  & = \mathbf{0} \quad \text{ on } \Omega \times (0, \infty) \label{eq:ivp} \\
    \mathbf{n} \times \ve & = \mathbf{0} \quad \text{ along } \partial\Omega \times (0, \infty) \label{eq:bc} \\
    \sigma \ve|_{t=0} & = \mathbf{q} \quad \text{ on } \Omega. \label{eq:source}
    \end{align}
\end{subequations}

\subsection{Finite element discretization}
We employ a finite element discretization based on Nédélec spaces on tetrahedral meshes.

\subsubsection{Variational formulation}
The variational formulation for Maxwell's equations seeks the solution of \eqref{eq:pde} in a suitable function space. 
An admissible space is the Sobolev space
\begin{equation}
    \mathbf H(\text{curl}; \Omega) = \left\{
    \vu \in L^2(\Omega)^3: \nabla \times \vu \in L^2(\Omega)^3
    \right\},
\end{equation}
where $L^2(\Omega)^3$ denotes the space of square-integrable vector-valued functions defined on $\Omega$.
By exploiting the boundary condition \eqref{eq:bc}, we can further restrict these fields to the subspace of  fields with vanishing tangential trace
\begin{equation}
    \mathcal V := \left\{
    \vu \in \mathbf{H}(\text{curl}; \Omega) : \mathbf n \times \vu = \mathbf{0} \text{ along } \partial\Omega
    \right\}.
\end{equation}
The variational formulation is obtained by  multiplying \eqref{eq:ivp}-\eqref{eq:source} with an arbitrary vector field $\Phi \in \mathcal V$.
After integration by parts, we arrive at the problem of seeking $\ve = \ve(\vr, t)$ such that
\begin{subequations}
\label{eq:weakform}
\begin{align}
    (\mu^{-1} \nabla \times \ve, \nabla \times \Phi) + (\sigma \partial_t \ve, \Phi) & = 0 \\
    (\sigma \ve|_{t=0}, \Phi) &= (\mathbf q, \Phi)
\end{align}
\end{subequations}
for all $\Phi \in \mathcal V$. The inner product on $L^2(\Omega)^3$ is denoted by $(\cdot, \cdot)$.

\subsubsection{Discretization in space}
We restrict the trial and test functions in the weak form \eqref{eq:weakform} to a finite-dimen\-sional subspace $\mathcal V^h \subset \mathcal V$ that consists of curl-conforming Nédélec elements on a tetrahedral mesh $\mathcal T_h$.
On each tetrahedron $T \in \mathcal T_h$, the functions in $\mathcal V^h$ consist of vector polynomials~$\vv$.
We obtain the discrete approximation of the solution of the variational problem \eqref{eq:weakform} by determining $\ve^h \subset \mathcal V^h$ such that for all $\Phi \in \mathcal V^h$ we have
\begin{subequations}
\begin{align}
        (\mu^{-1} \nabla \times \ve^h, \nabla \times \Phi) + (\sigma \partial_t \ve^h, \Phi) & = 0 \\
    (\sigma \ve^h|_{t=0}, \Phi) &= (\mathbf q, \Phi)
\end{align}
\end{subequations}

By expanding the discrete solution $\ve^h \in \mathcal V^h$ in a basis $\{\Phi_1, \dots, \Phi_N\}$ of $\mathcal V^h$, we arrive at the initial-value problem
\begin{subequations}
\label{eq:ode}
\begin{align}
    \mK \vu(t) + \mM \partial_t \vu(t) & = \mathbf{0}, \quad t \in (0, \infty) \\
    \mM \vu(0) & = \mathbf f.
\end{align}
\end{subequations}
The vector $\vu(t)$ contains $N$ coefficients of the finite-element approximation $\ve^h(t)$ with respect to the Nédélec basis at $t>0$. 
The curl-curl (stiffness) and mass matrices $\mK$ and $\mM$ as well as the initial values $\mathbf f$ are given in terms of the Nédélec basis by
\begin{subequations}
\begin{align}
    [\mK]_{i,j} & = (\mu^{-1} \nabla \times \Phi_j, \nabla \times \Phi_i) \\
    [\mM]_{i,j} & = (\sigma \Phi_j, \Phi_i) \\
    [\mathbf f]_{i,j} & = (\mathbf q, \Phi) \quad \text{for } i,j = 1,\dots, N.
\end{align}
\end{subequations}
The explicit solution of the semi-discretized, time-continuous problem \eqref{eq:ode} is given in  terms of the matrix exponential function acting on a vector,
\begin{equation}
    \vu(t) = \exp(-t \, \mM^{-1}\mK) (\mM^{-1}\vf) = \exp(-t \, \mM^{-1}\mK) \vb
\label{eq:expm}
\end{equation}
with $\vb = \mM^{-1}\vf$.

\section{Rational approximation}
\label{sec:numerical_method}

Following on from \eqref{eq:expm}, we need to derive a computationally viable method to compute approximations to $\mathbf{u}(t)$ at a number of, say, $K$ time channels $t_1,\ldots,t_K$. The use of rational approximants has been widely explored, in particular, when computed implicitly via rational Krylov projection methods; see, e.g., \citet{druskin2009solution,borner2015three,QG19}. In these approaches, the employed rational function is  computed \emph{implicitly} and hence difficult to work with in downstream tasks such as inversion. Furthermore, while approaches to parallelize rational Krylov methods exist, there are numerical stability constraints on the achievable level of parallelism.  

In Sections~\ref{sec:single} and \ref{sec:rkfit} below we present two approaches that are much more amenable to parallel implementation while delivering best or near-best approximation accuracies. To the best of our knowledge, both approaches have not been discussed in the context of TEM modeling. While we will argue that the second approach is superior, the first is still attractive as a baseline for comparison.

\subsection{Single-time best approximation} 
\label{sec:single}

A  direct approach to approximating the exponential is to employ a single rational best approximation $r_m(x) \approx \exp(-x)$ for $x\geq 0$. The theory and numerical computation of these approximants is a classical problem in approximation theory; see, e.g., \cite{cody1969chebyshev} and \cite[Chapter~25]{trefethen2019approximation}.  If we expand $r_m$ in partial fraction form
\[
r_m(x) = \alpha_0 + \sum_{i=1}^m \frac{\alpha_i}{x - \xi_i},
\]
then
\begin{subequations}
\begin{align}
r_m(t\, \mM^{-1}\mK)\vb 
&= \alpha_0 \vb + \sum_{i=1}^m \alpha_i   (t\,  \mM^{-1}\mK - \xi_i  \mI)^{-1}\vb  \\
&= \alpha_0 \vb + \sum_{i=1}^m \alpha_i  (t\,  \mK - \xi_i \mM)^{-1} \vf.
\label{eq:rba}
\end{align}
\end{subequations}
The evaluation of this formula for all time channels $t \in \{ t_1,\ldots,t_K\}$ requires $m\cdot K$ solves of complex shifted linear systems with the matrices $(t_j \mK - \xi_i \mM)$. This can be reduced to $K \lceil m/2\rceil$ solves by exploiting that all complex poles $\xi_i$ appear in conjugate pairs. For example, we will demonstrate below (Table~\ref{tab:errs}) that to approximate $\exp(-x)$ uniformly to $10^{-6}$ accuracy for all $x\in[0,+\infty)$, we require a degree of~$m=7$. For $K=31$ time channels and by exploiting the complex conjugacy of the poles, this amounts to factorizing and solving for $4\times 31=124$ complex shifted linear systems.

\subsection{RKFIT approach} 
\label{sec:rkfit}
There is an alternative approach that will only require about $\lceil \widehat m/2\rceil $ solves of complex shifted linear systems independent of $K$, where $\widehat m$ is usually a moderate multiple of $m$. The precise relation between $m$ and $\widehat m$ is complicated and depends on the ratio $t_{\max}/t_{\min}$, i.e., the relative length of the time interval one wants to cover. Table~\ref{tab:errs} shows some of these relations for selected interval ratios and target accuracies. For example, if $t_{\max}/t_{\min} = 10^3$ and the target  accuracy is $10^{-6}$, then $\widehat m=27$, and so in total only $14$ shifted linear systems need to be solved, \emph{independent} of $K$.

We follow the RKFIT approach described in~\cite[Section~6.2]{berljafa2017rkfit}.  In essence, the idea is to use a family of rational approximants \emph{with shared poles},  
\begin{equation}
r_{\widehat m}^{[j]}(x) =  \sum_{i=1}^{\widehat m} \frac{\alpha_i^{[j]}}{x - \xi_i} \approx \exp(-t_j x), \quad x\geq 0, \quad j=1,\ldots,K.
\label{eq:ratfam}
\end{equation}
RKFIT  (which stands for \emph{rational Krylov fitting}) computes such approximants  by using an implicit pole relocation technique to compute a sequence of rational Krylov spaces. The objective function to minimize approximately is of the form
\begin{equation}
\sum_{j=1}^K w_j \| r_{\widehat m}^{[j]}(\mS)\vv - \exp(t\mS) \vv \|_2^2 \to \min,
\label{eq:rkfit}
\end{equation}
where in our application the matrix $\mS$ is a diagonal surrogate matrix with ``dense'' eigenvalues on $[0,+\infty)$, $\vv$ is the vector of all ones, and the $w_j>0$ are weights that can be either all set to $1$, or chosen problem-dependent so that the weight varies across the time points $t_j$. 

We have dropped the absolute terms in \eqref{eq:ratfam} for convenience, leading to a family of subdiagonal rational approximants. More importantly, note that \emph{now the time dependence is entirely contained in the residue terms.} Hence the evaluation of 
\[
r_{\widehat m}^{[j]}(\mM^{-1}\mK)\vb = \sum_{i=1}^{\widehat m} \alpha_i^{[j]}  (\mK - \xi_i \mM)^{-1} \vf
\]
requires factorizations of the same system matrices for all $t_j$. Of course, these systems can be solved in parallel. Assuming that $\widehat m$ is even and all poles appear in complex conjugate pairs,\footnote{Some of the even-degree RKFIT approximants  we have computed have two real poles on the negative real axis. For such approximants,  two factorizations of real-shifted linear systems will be required, on top of the $\widehat m/2-1$ factorizations with complex shifts. Such approximants only appear for rather large degrees $\widehat m$ and the double real poles might be artifacts of computing in floating-point arithmetic. In any case, having two real poles is just as easy to implement as \eqref{eq:rmj} but would complicate the notation unnecessarily, hence we omit it.}
\begin{equation}
r_{\widehat m}^{[j]}(\mM^{-1}\mK)\vb =   2 \Re \left( \sum_{i=1}^{\widehat m/2} \alpha_i^{[j]}  (\mK - \xi_i \mM)^{-1} \vf \right).
\label{eq:rmj}
\end{equation}
Conveniently, by collecting the required $\widehat m/2$ solves in a  tall-skinny matrix
\[
\mS = \big [ (\mK - \xi_1 \mM)^{-1} \vf  , \ldots, (\mK - \xi_{\widehat m/2} \mM)^{-1} \vf\big ] 
\]
and the residues in 
\[
    \mR = \begin{bmatrix}
    \alpha_1^{[1]} & \cdots & \alpha_1^{[K]} \\
    \vdots & & \vdots \\
    \alpha_{\widehat m/2}^{[1]} & \cdots & \alpha_{\widehat m/2}^{[K]}
    \end{bmatrix},
\]
all $K$ time snapshots can be computed  from \eqref{eq:rmj} in one go as the columns of
the matrix-matrix product $\mS \mR$. 
This operation is rich in level-3 BLAS, which may add further performance benefits.

\subsection{Comparison with single-time best approximation}

We note that for $\mM$ symmetric positive definite and $\mK$ positive definite, we have
\[
\| \exp(-t_j \mM^{-1} \mK )\vb - r^{[j]}(\mM^{-1}\mK)\vb\|_{\mM} 
\leq 
\| \vb\|_{\mM} \max_{x\geq 0} \big| e^{-t x} - r^{[j]}(x) \big|
\]
with the vector $\mM$-norm $\| \vx\|_{\mM} := (\vx^T \mM \vx)^{1/2}$, 
for \emph{any} rational function $r^{[j]}$. A proof of this can be found, e.g., in  \cite[Thm.~3.2]{borner2015three}. Our aim should therefore be to minimize the scalar uniform approximation error on the right-hand side of this inequality across all time points of interest~$t_j$.

\begin{figure}
    \centering
    \includegraphics[width=0.95\linewidth]{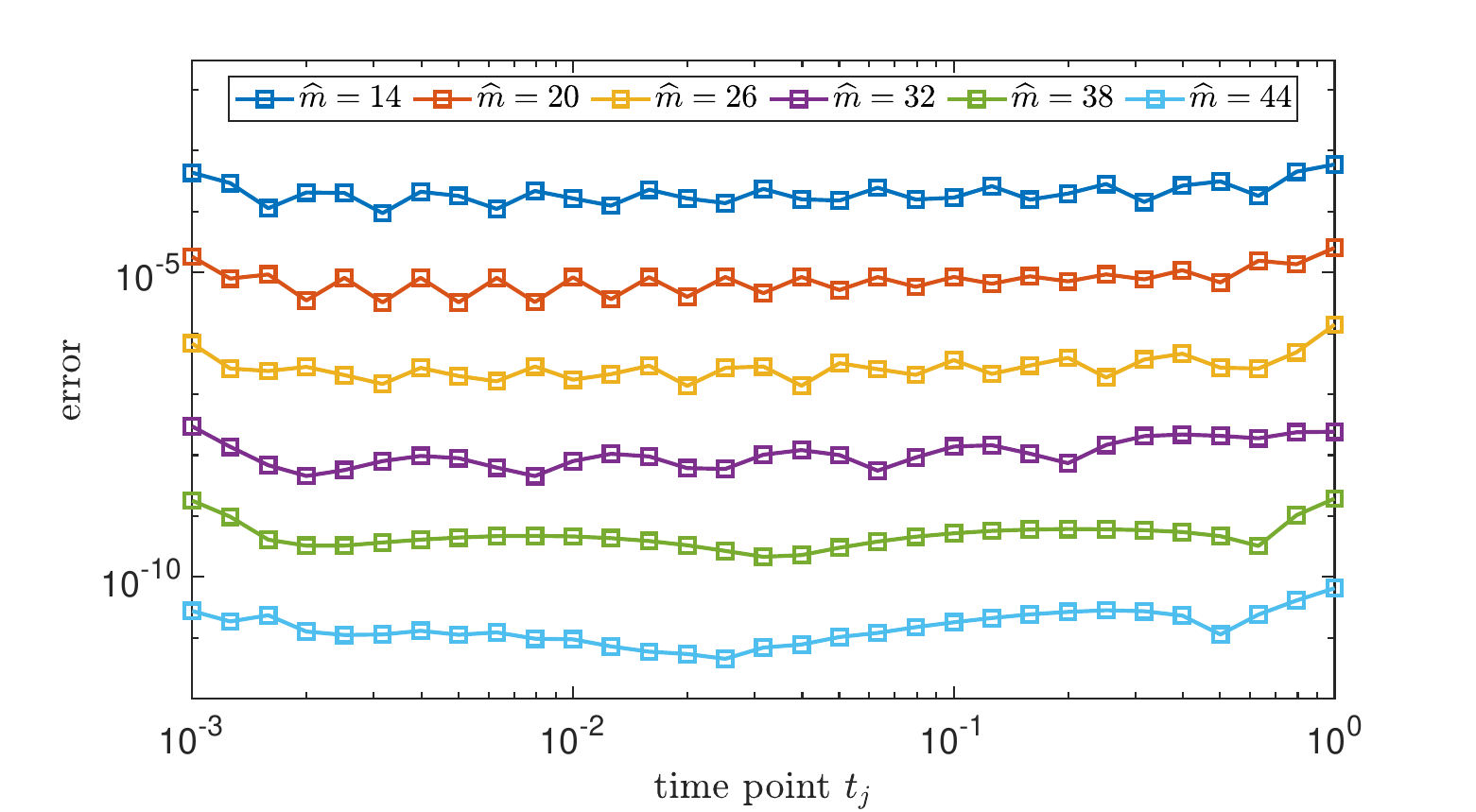}
    \includegraphics[width=0.95\linewidth]{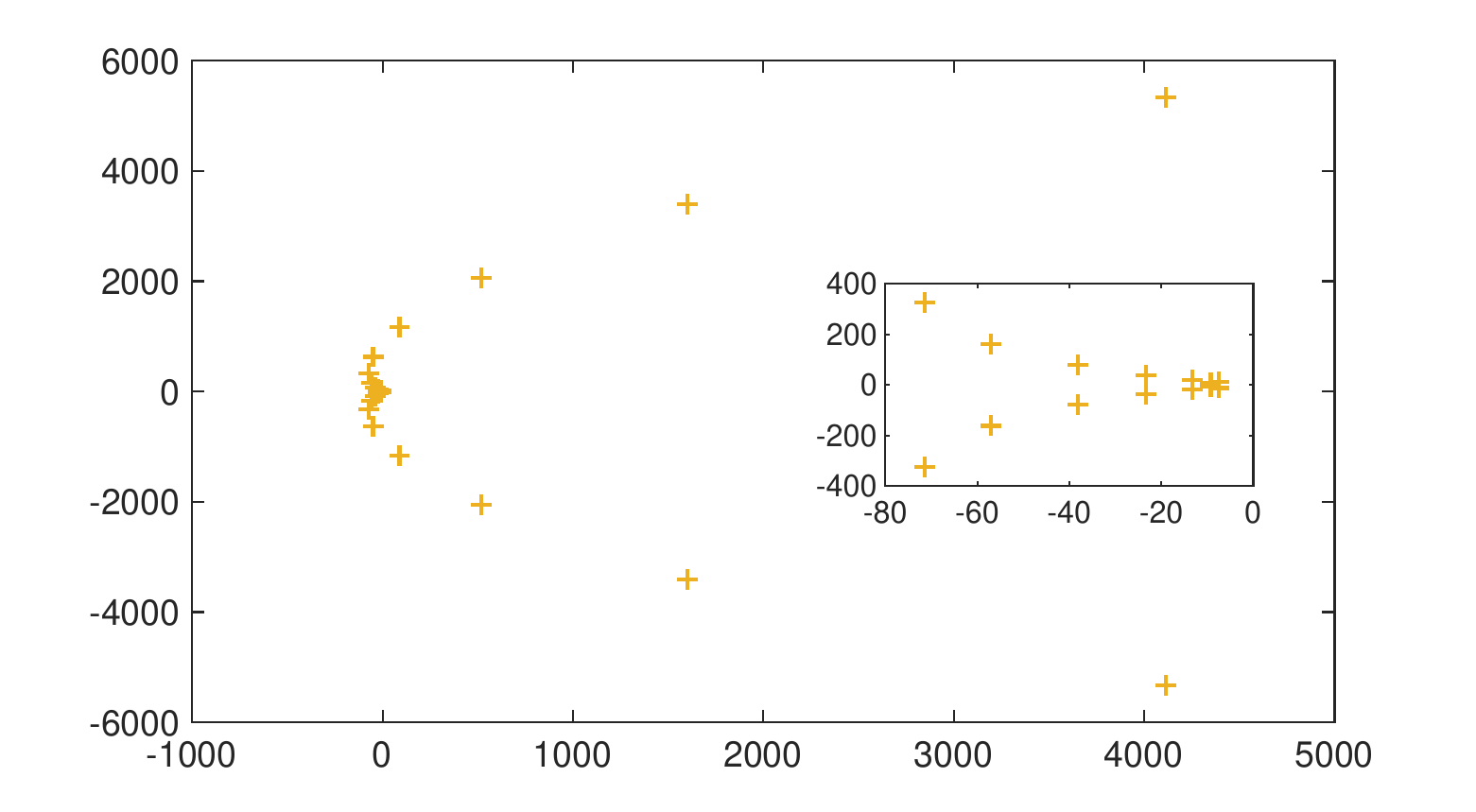}
    \caption{Top: Error of approximating $\exp(-tx)$ for $x\geq 0$ and for 31~logspaced time points $t_j \in [10^{-3}, 1]$. The degree of the shared-poles RKFIT approximants  is $\widehat m=14,20,\ldots,44$. Bottom: Poles of the RKFIT rational function family when the degree is  $\widehat m=28$, including a zoom close to the origin. In this case, the pair of  poles closest to the origin is $-7.574 \pm 12.45\mathrm{i}$.}
    \label{fig:err_plots1}
\end{figure}

\begin{figure}
    \centering
    \includegraphics[width=0.95\linewidth]{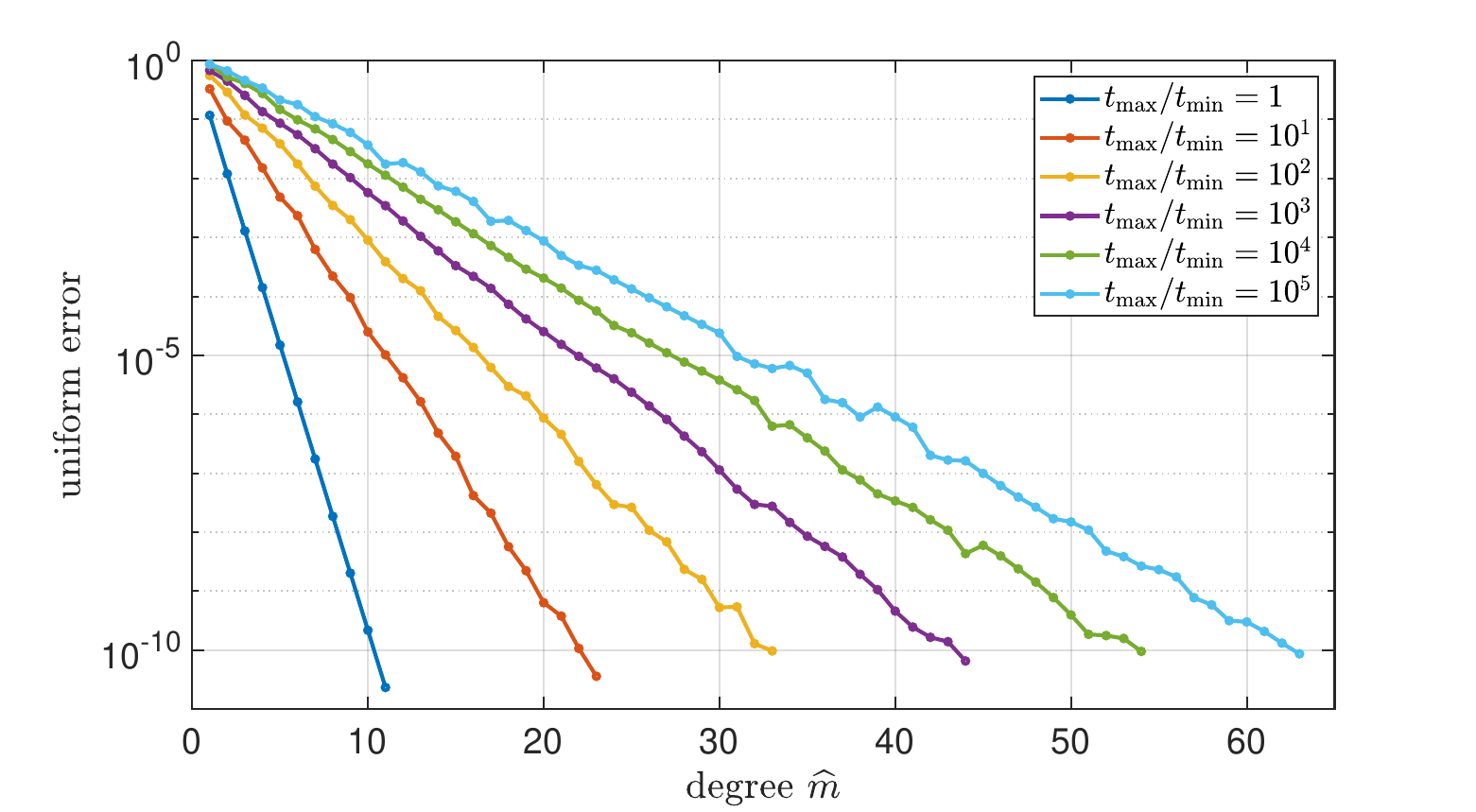}
    \caption{Error of the shared-poles RKFIT  approximants to $\exp(-tx)$ on $x\geq 0$ for varying degrees $\widehat m$ and varying time ratios $t_{\max}/t_{\min}$.}
    \label{fig:err_plots2}
\end{figure}

In Figures~\ref{fig:err_plots1} and \ref{fig:err_plots2} we illustrate the accuracy behavior of the RKFIT family of rational functions. In all these plots, the ``error'' corresponds to the quantity
\begin{equation}
\label{eq:error}
    \texttt{error} = \max_{x\geq 0} | \exp(-t_j x) - r^{[j]}(x)|
\end{equation}
for a specific time point $t_j$, whilst
\begin{equation}
\label{eq:uniform_error}
    \texttt{uniform error} = \max_{j=1,\ldots,K} w_j \max_{x\geq 0} | \exp(-t_j x) - r^{[j]}(x)|.
\end{equation}
The weights $w_j$ in the RKFIT objective function \eqref{eq:rkfit} are all set to $1$. In Figure~\ref{fig:err_plots1} (top), the time interval is fixed as $[t_{\min},t_{\max}] = [10^{-3},1]$ and we show the error for $K=31$ logspaced time channels on that interval for various degrees~$\widehat m$.
In Figure~\ref{fig:err_plots1} (bottom), we show the poles of the family of RKFIT approximants of degree~$28$. The plot confirms that these poles appear in 14 pairs of complex conjugate poles, and that they stay away from the spectral region $x\geq 0$, which would be particularly important if the corresponding linear systems with these shifts are solved with iterative solvers.

Table~\ref{tab:errs} and Figure~\ref{fig:err_plots2} explore the dependence of the uniform error as functions of the ratio $t_{\max}/t_{\min}$ and the degree~$\widehat m$.
We observe that, for a fixed time interval $[t_{\min},t_{\max}]$ the uniform error appears to decrease geometrically like $O(C^{-\widehat m})$ as $\widehat m$ increases ($C>1$). This seems plausible in view of what is known for the rational best approximation to $\exp(-x)$, i.e., for a single time point $t=1$. In this case, $\texttt{error}\sim 2H^{-m-1/2}$ with Halphen's constant $H = 9.289\ldots$; see \cite[Chapter~25]{trefethen2019approximation}. For the time-uniform RKFIT approximants, the geometric rate $C$ clearly depends on the ratio $t_{\max}/t_{\min}$. 

Assuming that, for a given target tolerance $\texttt{tol}$ we have $\texttt{err\_single} \sim H^{-m}$ and $\texttt{err\_rkfit} \sim C^{-\widehat m}$, we can estimate the required degrees of the single-time best rational and the time-uniform RKFIT  approximants  as
\[
m = -\log_H(\texttt{tol}), \quad \widehat m = -\log_C(\texttt{tol}),
\]
respectively. 
The quotient $\widehat m/m$ is independent of $\texttt{tol}$. When evaluating  approximations to the matrix exponential for $K$ time points and using direct solvers, we need $K \lceil m/2 \rceil$ matrix factorizations with the single-time best rational approximant, while our new time-uniform RKFIT family only requires $\lceil \widehat m/2\rceil$ matrix factorizations independent of $K$. Hence, the critical value for which the RKFIT family will become more efficient to evaluate is
\[
K^* = \lceil \widehat m/2 \rceil \big/ \lceil m/2 \rceil.
\]
At the bottom of Table~\ref{tab:errs}, we list computed averages of this critical value for each ratio $t_{\max}/t_{\min}$. For example, if $t_{\max}/t_{\min} = 10^3$, then $K^* \approx 4.1$ and so evaluating our family of time-uniform approximants at $K \geq 5$ time points will require fewer matrix factorizations than evaluating $K$ single-time rational approximants.

\begin{table}
    \centering
    \begin{tabular}{l  c  c  c  c  c  c}\toprule
           accuracy  & \multicolumn{6}{c}{$t_{\max}/t_{\min}$}\\
              \cmidrule(lr){2-7}
              & 1   & $10^1$ & $10^2$ & $10^3$ & $10^4$ & $10^5$ \\\midrule
         $10^{-2}$ & 2   & 5      & 7      & 10     & 12     & 14 \\
         $10^{-4}$ & 4   & 9      & 14     & 18     & 22     & 26 \\
         $10^{-6}$ & 7   & 14     & 20     & 27     & 33     & 38 \\
         $10^{-8}$ & 9   & 18     & 27     & 35     & 44     & 52 \\ 
        $10^{-10}$ & 11  & 23     & 33     & 44     & 54     & 63 \\\midrule
        $\texttt{avg}(K^*)$ & $-$ &  2.2 & 3.1 &  4.1 & 4.9 &  5.8  \\\bottomrule  
    \end{tabular}
    \caption{Degree $m$ required to achieve the target accuracy shown in the left column over a time interval $[t_{\min},t_{\max}]$ with $t_{\max}/t_{\min}$ shown in the first row. The degrees for $t_{\max}/t_{\min}=1$ are those for the best approximation to the exponential function of type $[m/m]$ known by \protect\cite{cody1969chebyshev}, while the other columns are approximants of type $[(\widehat m-1)/\widehat m]$ computed using RKFIT.}
    \label{tab:errs}
\end{table}

\section{Implications for TEM inversion}
\label{sec:inversion}

The inverse TEM problem consists of recovering the electrical conductivity $\sigma$ from measurements of the transient response.
For the matrix exponential formulation \eqref{eq:expm}, this means recovering the parameters $\mathbf m = [ m_1; m_2 ; \ldots; m_P ]$ in 
$$
\vu(t, \vm) = \exp(-t\, \mM(\vm) \mK)(\mM(\vm)^{-1}\vf)
$$
from measurements $\mQ \vu(t_j)$, where $\mQ\in\mathbb{R}^{M\times N}$ with $M\leq N$. We will assume that 
\[
\mM(\vm) = e^{m_1}\mM_1 + e^{m_2}\mM_2 + \cdots + e^{m_P} \mM_P
\]
with some fixed matrices $\mM_1,\mM_2,\ldots,\mM_P$. 
Differentiating $\vu(t,\vm)$ with respect to $\vm$ is inconvenient for it involves the product of two $\vm$-dependent terms. The problem simplifies considerably if, instead, we consider the rational approximation
\begin{eqnarray}
\vv(t_j, \vm) &:=&  2 \Re \left(  \sum_{i=1}^{\widehat m/2} \alpha_i^{[j]} \mQ (\mK - \xi_i \mM(\vm))^{-1} \vf \right) \nonumber\\
&\,=& 
2 \Re \left(  \sum_{i=1}^{\widehat m/2} \alpha_i^{[j]} \mQ \mA_i(\vm)^{-1} \vf \right),
\label{eq:rmjQ}
\end{eqnarray}
where we have denoted $\mA_i(\vm) :=  \mK - \xi_i \mM(\vm)$. We will use the identity
\[
\partial_\vm \mA_i(\vm)^{-1} = -\mA_i(\vm)^{-1} \left[\partial_\vm \mA_i(\vm) \right] \mA_i(\vm)^{-1},
\]
where the term $[\,\cdot\,]$ is an $N\times N \times P$ tensor of partial derivates and its product with matrices is page-wise. The pages of $[\,\cdot\,]$ are given as
\[
\left[\partial_\vm \mA_i(\vm) \right](\,:\,,\,:\,,k) 
=
- \xi_i \, e^{m_k} \mM_k, \quad k = 1,\ldots,P.
\]
The Jacobian of $\vv$ at $t=t_j$ with respect to $\vm$ is
\[
\mJ_j := \partial_{\vm} \vv(t_j,\vm) 
= 
2 \Re \left( \sum_{i=1}^{\widehat m/2} \alpha_i^{[j]} \xi_i \mQ \mA_i(\vm)^{-1} \big[e^{m_1}\mM_1 \vg_i,\ldots, e^{m_P}\mM_P \vg_i \big]\right),
\]
where $\vg_i := \mA_i(\vm)^{-1} \vf$ needs to be computed only once for each $i$. The required factorization of $\mA_i(\vm)$ is normally required anyway for the forward evaluation of \eqref{eq:rmjQ}, and it can also be reused for the block solve with the matrix appearing in $[\,\cdot\,]$. Hence, with direct solvers, the evaluation of the forward solution as well as its full Jacobian requires only $\widehat m/2$ matrix factorizations, $\widehat m/2$ solves with right-hand side vector $\vf$, and $\widehat m/2$  solves with a block matrix having $P$ columns. 

Further note that the only time dependence when evaluating $\partial_{\vm} \vv(t_j,\vm)$  is in the residues $\alpha_i^{[j]}$ which only enter as multiplicative factors. As a consequence, we can get the full Jacobian across all time points $t_1,\ldots,t_K$ as
\[
\mJ:= \begin{bmatrix}
    \mJ_1 \\
    \mJ_2 \\
    \vdots \\
    \mJ_K
\end{bmatrix}
=
2 \Re \left( \sum_{i=1}^{\widehat m/2} \begin{bmatrix}
\alpha_i^{[1]} \\
\alpha_i^{[2]} \\
\vdots \\
\alpha_i^{[K]} 
\end{bmatrix} 
\otimes 
\left(
\xi_i \mQ \mA_i(\vm)^{-1} \big[e^{m_1}\mM_1 \vg_i,\ldots, e^{m_P}\mM_P \vg_i \big]\right)
\right).
\]
This is an $MK \times P$ matrix whose evaluation requires only $\widehat m/2$ matrix factorizations, $\widehat m/2$  solves with right-hand side vector $\vf$, and $\widehat m/2$  solves with a block right-hand side having $P$ columns.  

We can readily use the Jacobian within a Gauss--Newton solver \cite{kelley2003solving}. One Gauss--Newton step is   given as
\[
\vm_{\ell+1}  = \vm_\ell
- 
\alpha\cdot \mJ^\dagger (\vv(:,\vm_\ell) - \vv(:,\vm_0)),
\]
where $\vv(:,\vm_0),\vv(:,\vm_\ell)$ are vectors of length $MK$ collecting the time observations for the ground truth $\vm_0$ (the sensor measurements) and the current parameter guess $\vm_\ell$, respectively, and $\mJ^\dagger$ stands for the Moore--Penrose pseudoinverse (i.e., a least squares solution) which can be obtained via economic QR decomposition of~$\mJ$. The time step parameter $\alpha>0$ can be fixed or chosen via a line search method.

\subparagraph{Using iterative solvers.} In most cases, it may be more suitable to not form the full $MK\times P$ Jacobian~$\mJ$ as we did above, especially if the systems $\mA_i(\vm)$ are sparse and so large that direct solvers are not applicable and/or the number of parameters~$P$ is very large. 
In this case, one can apply, for example, CGLS which will merely require the actions $\mJ \vx$ and $\mJ^T \vy$ on vectors \cite{bjorck}. The implementation of this still only requires $\widehat m/2$ matrix factorizations per Gauss--Newton iteration, and each action $\mJ \vx$ (or $\mJ^T \vy$) requires $2\times \widehat m/2$ solves with $\mA_i$ (or~$\mA_i^T$). 

In summary, we conclude that using our new time-uniform RKFIT approximants significantly reduces the cost of both forward TEM simulation and inversion, and the cost for solving linear systems of equations (usually the computational bottleneck) is independent of the number of time channels $K$. Further, as the rational approximants are uniformly good over the whole real half axis $x\geq 0$, the number of required solves is mesh independent and it will also not depend on conductivity contrasts.

This concludes the basic derivation of an  inverse solver. We follow with two brief comments on possible extensions.

\subparagraph{Regularization.}
In some cases, one wishes to impose that $\vm$ remains close to some reference model $\vm_r$, or that $\vm$ exhibits some smoothness. This can be achieved by extending $\mJ$ by a $P\times P$ matrix $\sqrt{\lambda}\, \mI$ and extending the vectors $\vv(:,\vm_0)$ and $\vv(:,\vm_\ell)$ by $\sqrt{\lambda}\vm_r$ and  $\sqrt{\lambda} \vm_\ell$, respectively.

\subparagraph{Using multiple sources.} The above derivation straightforwardly generalizes to using a block source $\vf\in\mathbb{R}^{N\times S}$. In this case, the full Jacobian $\mJ$ becomes an $MKS\times P$ matrix.

\section{Numerical experiments}
\label{sec:experiments}

\begin{figure}
    \centering
    \includegraphics[width=0.65\linewidth]{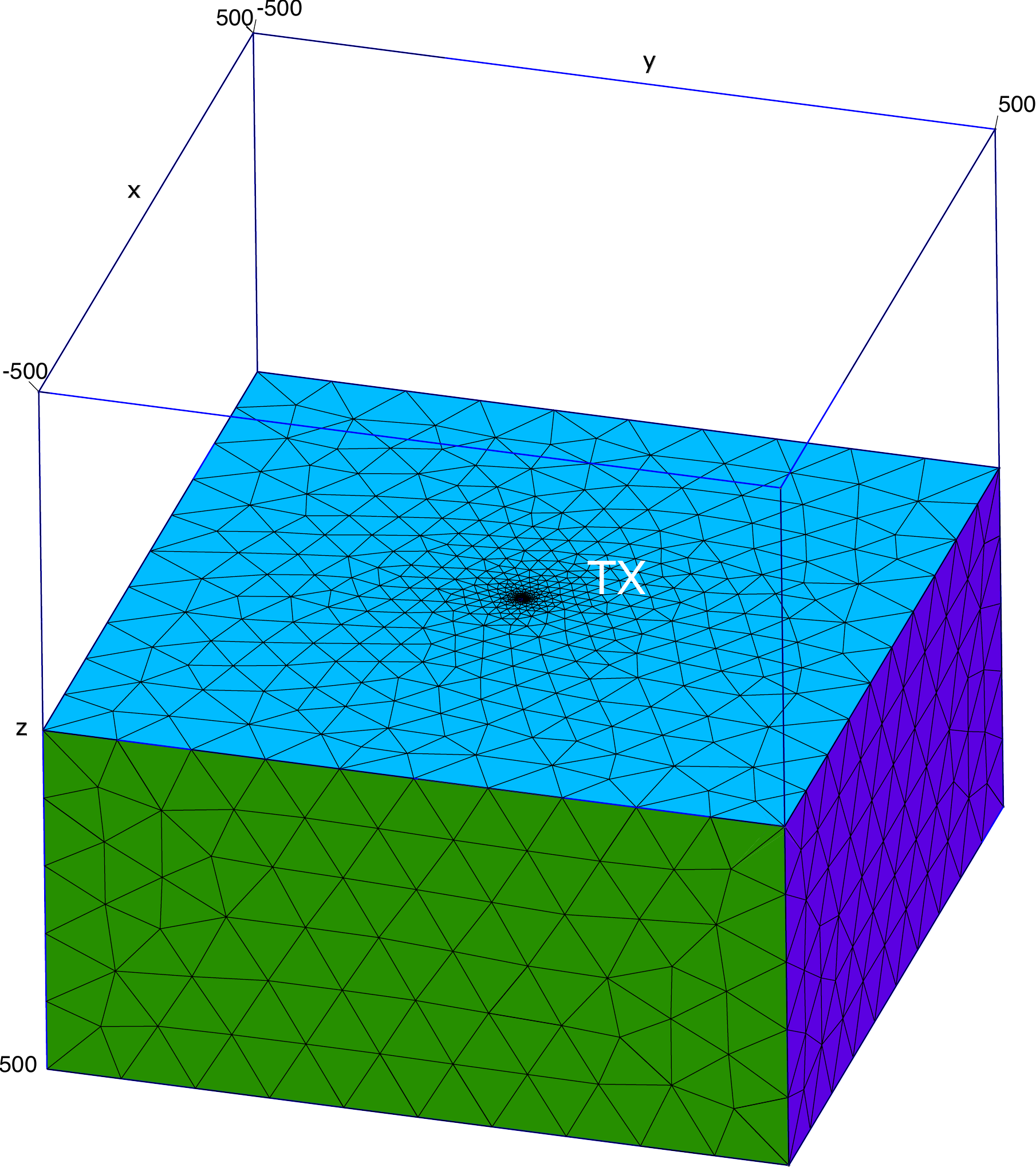}
    \caption{Partial view of the tetrahedral mesh used for the numerical experiments. At at $z=0$, a square loop of side length 5 m (TX) is located at the top of the conducting halfspace overlain by the air layer (not shown).}
    \label{fig:mesh}
\end{figure}

To evaluate the performance of the proposed numerical method, we consider a model comprising a uniform conducting half-space with electrical conductivity $10^{-1}$~S/m, overlain by a resistive air half-space with conductivity $10^{-8}$~S/m.
A square horizontal transmitter loop of side length $L = 5$~m is positioned at the interface between the air and the conductive half-space (Fig.~\ref{fig:mesh}).
We compute the transient response, i.e., the time derivative of the vertical magnetic flux density, $\partial_t b_z(t)$, at the center of the loop for 31 logarithmically spaced time points within the interval $t \in [10^{-6}, 10^{-3}]$~s, following current shut-off at $t = 0$~s. 
The spatial discretization is performed using first-order Nédélec edge elements, resulting in a system with stiffness and mass matrices $\mathbf{K}, \mathbf{M} \in \mathbb{R}^{N \times N}$, and a source vector $\mathbf{f} \in \mathbb{R}^N$, with $N = 81,\!174$ degrees of freedom.

Note that we deliberately restrict our experiments to a simple model.
While more geologically realistic models could be considered, the focus of this study is not on the interpretation of complex subsurface responses, but rather on the numerical behavior of the forward solver.
Despite its simplicity, the chosen setup is sufficiently representative from a numerical standpoint: the spatial discretization introduces a wide range of element sizes, while the strong conductivity contrast between the air and the subsurface (seven orders of magnitude) ensures that the system exhibits an eigenvalue spectrum comparable to those encountered in practical 3-D scenarios.
As such, the test problem captures the essential numerical challenges relevant to more complex settings, allowing us to meaningfully assess the accuracy and efficiency of the proposed method.

Although the rational approximation directly yields an estimate of the electric field $\mathbf{e}(\mathbf{r}, t)$, the quantity of primary interest in most TEM applications is the time derivative of the magnetic flux density, $\partial_t \mathbf{b}(\mathbf{r}, t)$.
To extract this from the numerical solution, we implement an observation operator $\mathbf{Q} \in \mathbb{R}^{N}$ that evaluates a linear combination of the degrees of freedom in $\mathbf{u}(\mathbf{r}, t)$, corresponding to the curl of the Nédélec basis functions within the tetrahedral element containing the point $\mathbf{r}_0$ at which the response is desired.
Specifically, we approximate the vertical component of the curl of the electric field via
\[
-\partial_t b_z(\vr_0, t) =  \mathbf e_z \cdot (\nabla \times \mathbf e(\vr_0, t)) \approx \mathbf Q^\top(\vr_0) \vu(\vr, t),
\]
where $\mathbf{e}_z$ denotes the unit vector in the vertical direction. 

All numerical experiments were performed on a shared-memory computer with a non-uniform memory access (NUMA) architecture, equipped with four Intel® Xeon® Gold 5418Y CPUs (each with 24 cores, 96 cores in total) and 2 TB of RAM. 
The forward solver is implemented in the Julia programming language \citep{bezanson2017} using the \textit{Gridap.jl} finite element library \citep{Badia2020,Verdugo2022}, which provides a flexible and efficient framework for high-level implementation of finite element methods.
The computational mesh was generated using the Gmsh library \citep{geuzaine2009}, an open-source tool for finite element mesh generation and preprocessing. 
As a reference solution, we employed a one-dimensional model computed using the open-source library \textit{empymod} \citep{werthmuller2017}, which provides semi-analytical solutions for electromagnetic problems in layered media.


\subsection{Single-time best approximation}

\begin{figure}
    \centering
    \includegraphics[width=0.95\linewidth]{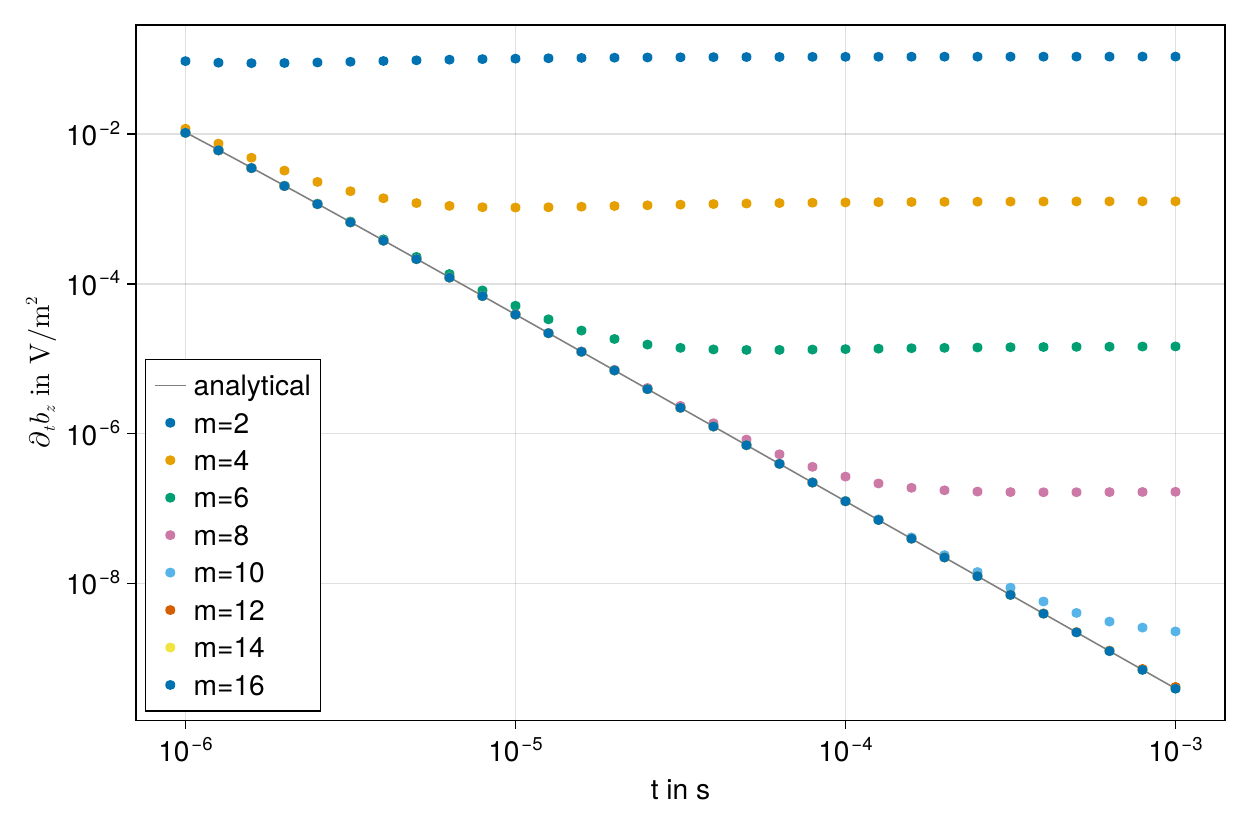}
    \caption{Single-time best rational approximation: Transient $\partial_t b_z(t)$ evaluated at the center of the transmitter loop for varying degrees $m=2,4,\dots,16$ and $K=31$ log-equidistant times compared to the analytical solution for the uniform half-space (solid line).}
    \label{fig:transientsRBA}
\end{figure}

The first numerical experiment illustrates the effect of the approximation degree $m$ on the accuracy of the single-time rational best approximation used in Section~\ref{sec:single}. 
As shown in Fig.~\ref{fig:transientsRBA}, the computed transients converge progressively towards the analytical reference solution as $m$ increases from $2,4,\ldots,10$.
Notably, the error exhibits geometric decay with increasing $m$, particularly at late times.
For degrees $m \geq 12$, the numerically computed transient is visually indistinguishable from the analytical solution. 
When looking at the relative error in Fig.~\ref{fig:strba_rkfit} (bottom), a degree of $m\geq 14$ is needed to capture the transient over the whole time interval of interest to high accuracy, in particular at late times. 

{We emphasize that for $m = 14$, the total number of linear system solves is $K \times m/2$, which amounts to $31 \times 7 = 217$ complex solves for $K = 31$, accounting for the fact that the poles $\xi_i$ occur in complex conjugate pairs.}

\subsection{RKFIT rational approximation}

We now compare the RKFIT-based rational approximants introduced in Section~\ref{sec:rkfit} against the optimal single-time rational approximation, where ''optimal'' refers to the lowest degree $m = 14$ for which the approximation error is sufficiently small across all evaluation times.

\begin{figure}
    \centering
    \includegraphics[width=0.95\linewidth]{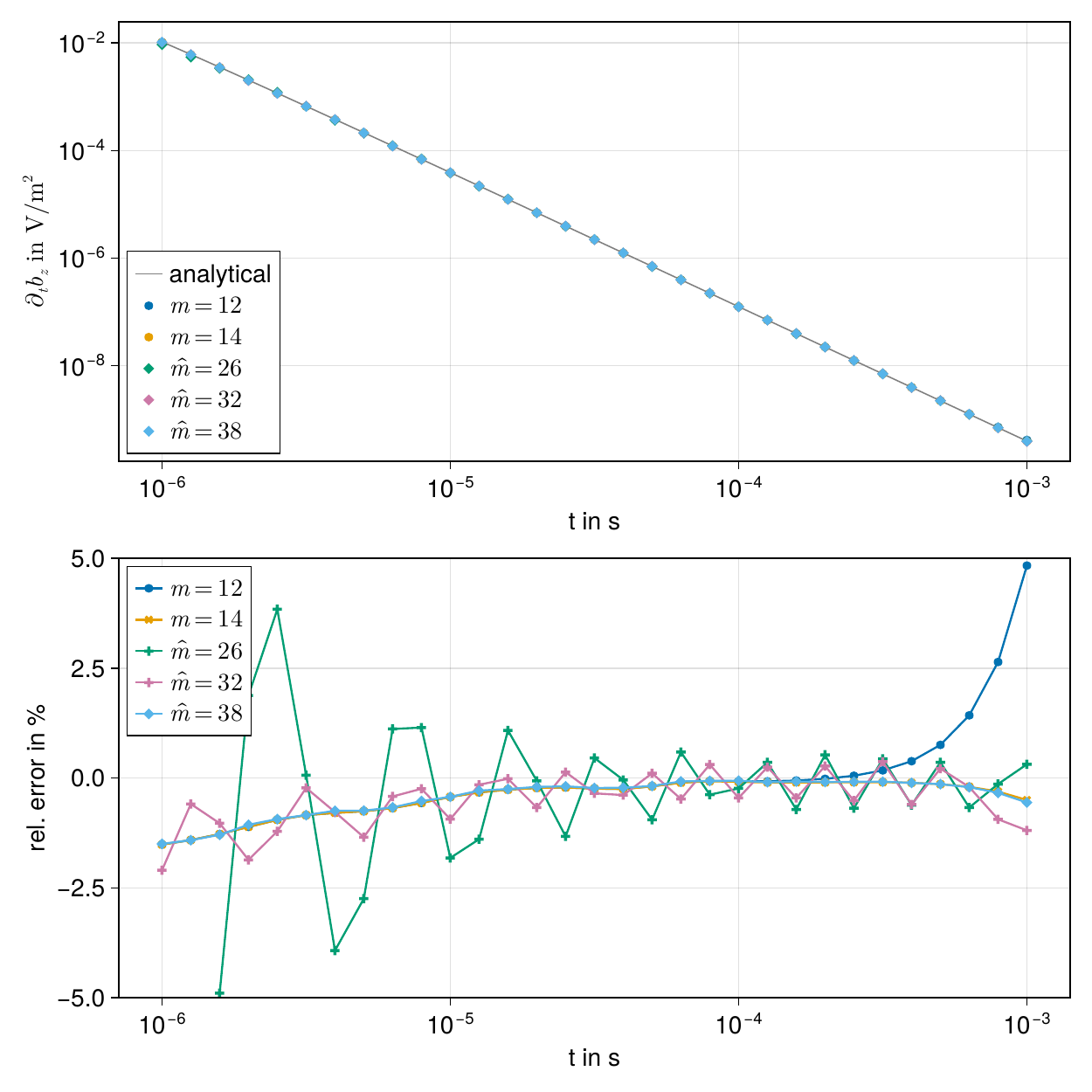}
    \caption{Accuracy comparison of the single-time best rational approximation for $m=12, 14$ and the RKFIT rational approximation for $\widehat m = 26, 32, 38$ with respect to the analytical solution evaluated at the center of the transmitter loop. Relative errors (in percent) are shown in the bottom panel.}
    \label{fig:strba_rkfit}
\end{figure}

Figure~\ref{fig:strba_rkfit} shows the  transients computed using the single-time rational approximation of degree $m = 12$ and $14$, and RKFIT-based approximants with degrees $\widehat{m} = 26$, $32$, and $38$.
The RKFIT approximants lead to computed transients that are comparable in accuracy to that of the single-time best rational approximants; see Fig.~\ref{fig:strba_rkfit} (top). When looking at the relative error in Fig.~\ref{fig:strba_rkfit} (bottom), RKFIT approximants appear slightly less smooth due to the uniform-in-time nature of the approximation. 
Nevertheless, an RKFIT approximant of degree $\widehat{m} = 38$ already achieves a level of accuracy on par with the single-time best rational approximation for $m = 14$.

We emphasize that for $\widehat m = 38$, the total number of linear system solves required with the RKFIT approach is just $\widehat m/2=19$, independently of $K$. When compared to the 217 complex solves from above, this is an improvement by a factor of more than~11.

We also investigated the influence of the weights $w_j$ in the RKFIT objective function \eqref{eq:rkfit} on the accuracy of the numerical solution at late times. The weights can be designed to offset the increased \emph{relative} numerical error that arises due to the decay of the fields at late times. 
Specifically, we conducted a series of experiments using weights $w_j \sim (t_j^{3/2}, t_j^{5/2})$, selected based on the asymptotic decay characteristics of the electromagnetic fields. 

Our weight choices are physically motivated. For a horizontal circular current loop of radius~$a$, centered on the $z$-axis of a cylindrical coordinate system, the horizontal electric field and the vertical magnetic field measured at $z=0$ and $r>0$ both decay as $t^{-3/2}$ in the late-time regime. In contrast, the time derivative of the vertical magnetic flux density, $\partial_t b_z$, decays more rapidly, following a $t^{-5/2}$ asymptotic behavior \citep{nabighian1988electromagnetic}. 
Our experience is that the lowest relative approximation errors for $\partial_t b_z(\mathbf{r}_0, t)$ at late times are indeed obtained for $w_j \sim t_j^{5/2}$.

\subsection{Parallel performance}

\begin{figure}
    \centering
    \includegraphics[width=0.95\linewidth]{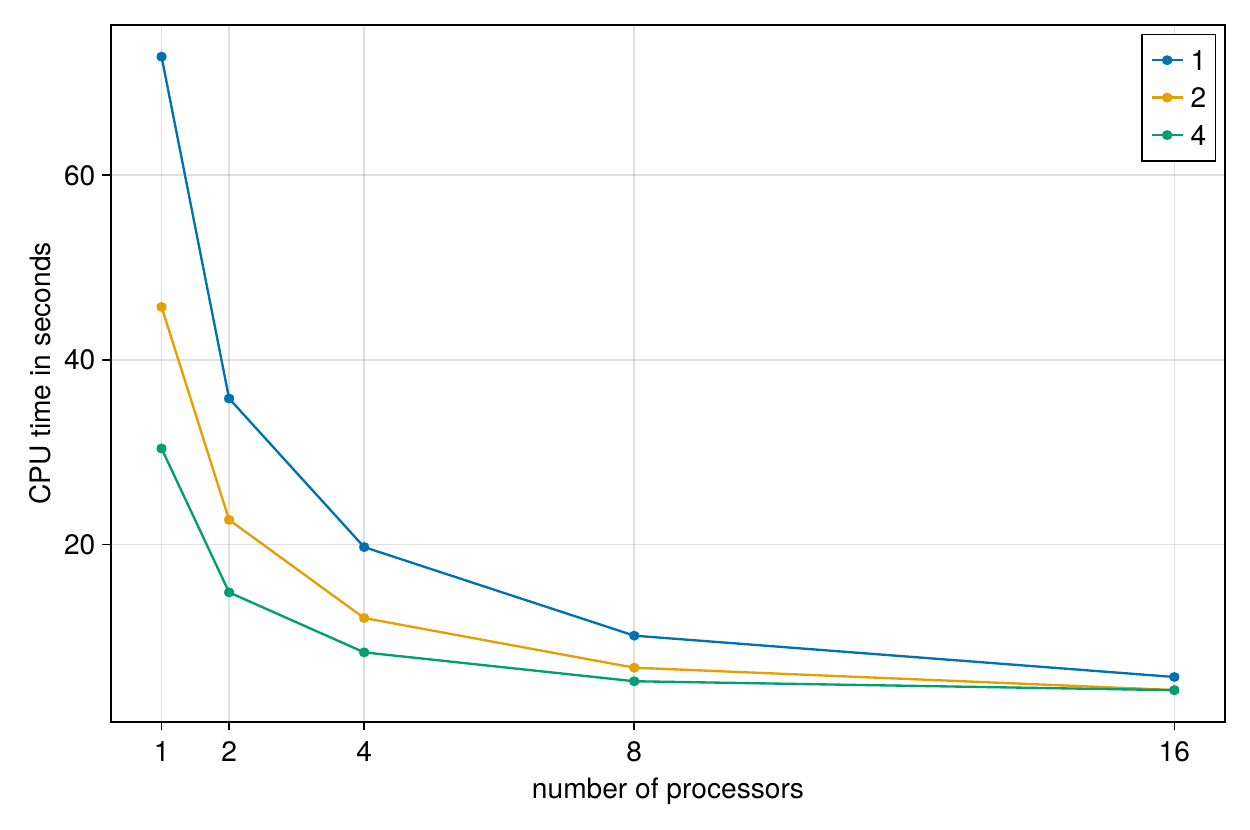}
    \includegraphics[width=0.95\linewidth]{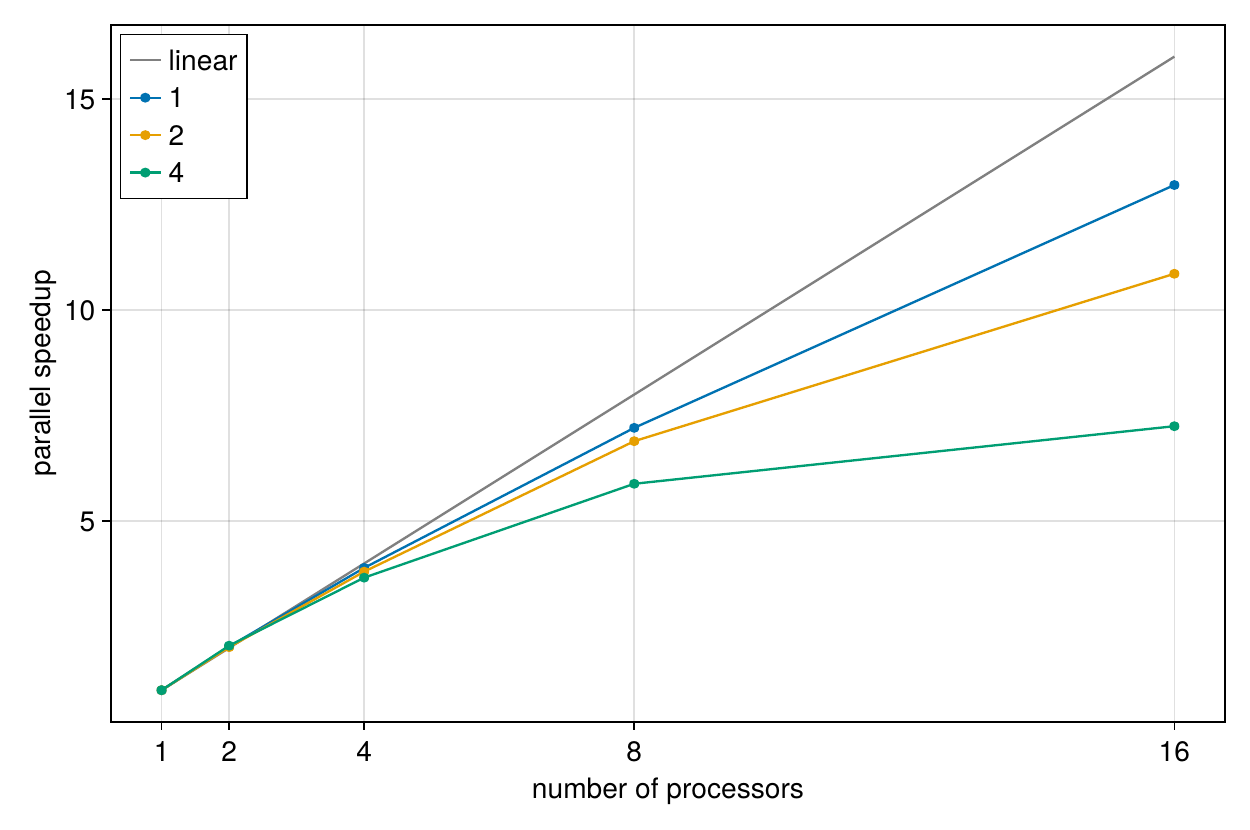}
    \caption{RKFIT rational approximation: CPU time (top) and parallel speedup with respect to the number of processors and the number of threads $n=1,2,4$ used by the PARDISO direct solver for degree $\widehat m=32$ and $K=31$.} 
    \label{fig:speedup}
\end{figure}

We investigate the parallel performance of the RKFIT rational approximants for transient electromagnetic modeling by distributing the linear systems across multiple processors for $\widehat{m}/2=16$. 
Numerical experiments were conducted using 1, 2, 4, 8, and 16 processors, and the corresponding wall-clock times are reported in Fig.~\ref{fig:speedup} (top). 
Additionally, we assess the impact of shared-memory parallelism by utilizing the multithreading capabilities of the Intel MKL PARDISO linear solver \citep{schenk2004solving} with 1, 2, and 4 threads.

Fig.~\ref{fig:speedup} (bottom) illustrates the observed speedup relative to the baseline serial execution and compares it to the ideal linear speedup. The speedup has been computed as
\[
S = \frac{T_1}{T_p},
\]
where $T_1$ denotes the execution time of the serial implementation, and $T_p$ is the execution time when using $p$ processors.
In all cases, we observe a clear reduction in computational time with increasing parallel resources. 
This demonstrates the scalability and efficiency of the approach in both distributed and shared memory settings.

Although our proposed algorithm is, in principle, embarrassingly parallel and expected to exhibit near-ideal scalability (in particular, when PARDISO is used in single-threaded mode), our experiments do not quite achieve optimal speedup.
We attribute this discrepancy to various factors, including the specific implementation of Julia’s \texttt{pmap} parallel map which introduces overhead in task scheduling and communication that limits parallel efficiency, particularly at higher processor counts.

Further deviation from the ideal speedup can be observed when the number of PARDISO threads is high (e.g., 4) and the number of processors is large (e.g., 16).
This shortfall is attributable to the suboptimal parallel scaling of PARDISO, not the RKFIT family of approximants whose partial fraction terms can be evaluated with perfect parallel efficiency. This problem may be attributed to memory bandwidth limitations or cache contention arising from oversubscription of shared hardware resources.
When multiple threads simultaneously access and manipulate large sparse data structures, such as those involved in the factorization and substitution steps of the linear solver, cache coherence overhead and reduced memory locality can significantly impact performance.

\begin{figure}
    \centering
    \includegraphics[width=0.95\linewidth]{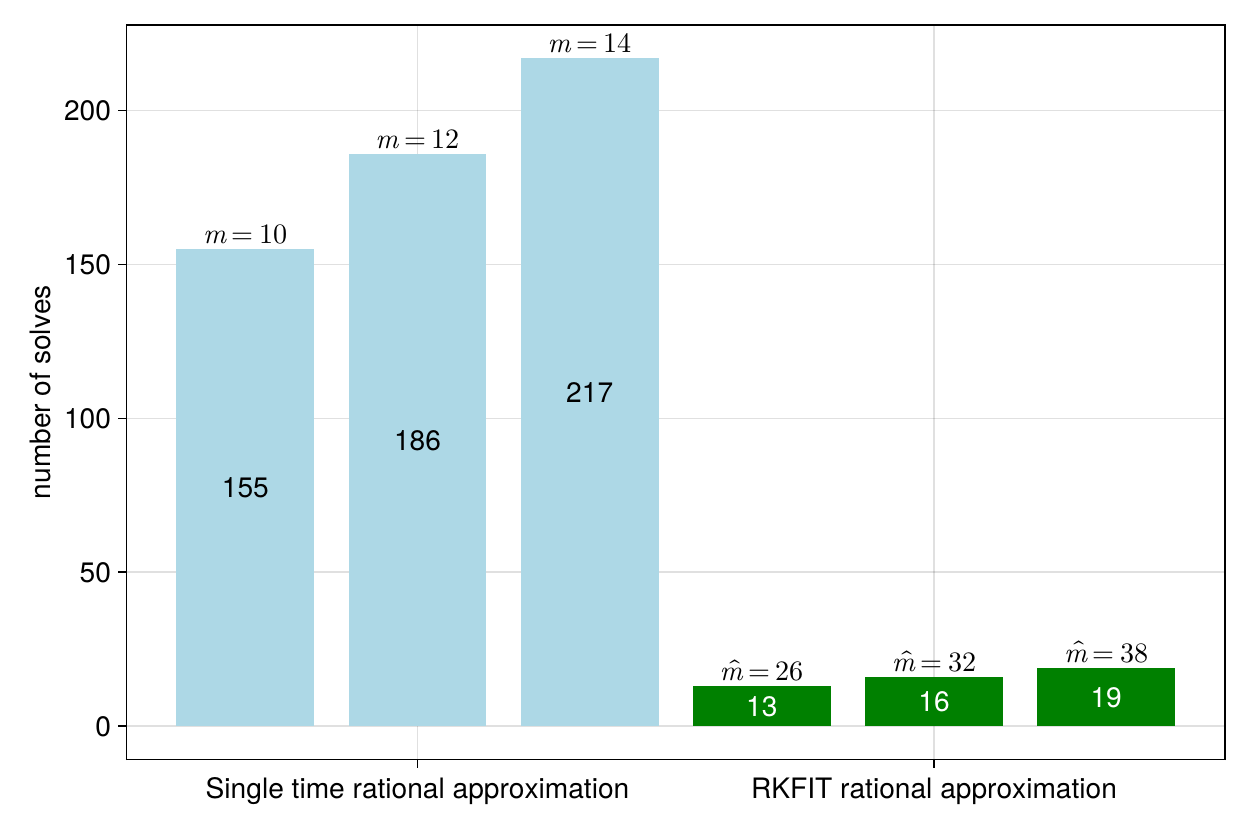}
    \caption{Number of complex linear system solves for varying degrees $m$ (single time rational approximation) and $\widehat m$ (RKFIT rational approximation), resp., for $K=31$ log-equidistant times.}
    \label{fig:solves}
\end{figure}

Finally, Fig.~\ref{fig:solves} presents a bar chart comparing the number of linear system solves required for the single-time best rational approximants of degrees $m = 10$, $12$, and $14$, against those required by the RKFIT rational approximants of degrees $\widehat{m} = 26$, $32$, and $38$, exploiting that poles appear in complex conjugate pairs.
In the experiments conducted with $m = 14$ and $\widehat{m} = 38$ we observe that RKFIT reduces the computational effort in terms of required system solves by more than a factor of 11 while maintaining a comparable level of accuracy. 
The results clearly show that RKFIT achieves comparable accuracy with significantly fewer system solves, highlighting its computational efficiency relative to the single-time best rational approximation approach.

\section{Conclusion}
\label{sec:conclusion}

We have introduced a highly parallel approach for the fast forward modeling of TEM fields. The approach exploits a family of rational approximants with time-independent shared poles computed by the RKFIT method. Our numerical experiments indicate that this approach is far superior over evaluating rational best approximants for each time point, even when the number of time channels is small (such as $K=5$). In a more practically relevant scenario where $K=31$, we observe parallel speedups of more than a factor of $10$. 

The applicability of our approach is not limited to TEM modeling. In principle, any linear parabolic evolution problem that needs to be solved over a given time interval could make use of the new rational function family and its parallel evaluation. In future work, we would also like to explore in detail how our approach performs in practice within an inversion loop and in combination with iterative linear system solvers.

\section*{Acknowledgments}
SG is supported the Royal Society (RS Industry Fellowship IF/R1/231032) and the Engineering and Physical Sciences Research Council (EPSRC grant EP/Z533786/1).

\bibliographystyle{gji}
\bibliography{refs}

\end{document}